\documentclass[12pt]{article}
\usepackage{amsmath, latexsym, amsfonts, amssymb, amsthm, amscd}
\usepackage{graphics}
\newtheorem{theorem}{Theorem}

\newtheorem{proposition}{Proposition}
\newtheorem{lemma}{Lemma}
\newtheorem{definition}{Definition}

\begin{document}

\title{Levy preservation and associated properties for $f$-divergence minimal equivalent martingale measures}

\maketitle

\begin{center}
{\large S. Cawston}\footnote{$^{,2}$ LAREMA, D\'epartement de
Math\'ematiques, Universit\'e d'Angers, 2, Bd Lavoisier - 49045,
\\\hspace*{.4in}{\sc Angers Cedex 01.}

\hspace*{.05in}$^1$E-mail: suzanne.cawston@univ-angers.fr$\;\;\;$
$^2$E-mail: lioudmila.vostrikova@univ-angers.fr}{\large
 and  L. Vostrikova$^2$}
\end{center}
\vspace{0.2in}

\abstract{We study such important properties of $f$-divergence minimal martingale measure  as Levy preservation property, scaling property, invariance in time property for exponential Levy models. We give some useful decomposition for $f$-divergence minimal martingale measures and we answer on the question which form should have $f$ to ensure mentioned properties. We show that $f$ is not necessarily common  $f$-divergence. For common $f$-divergences, i.e. functions verifying $f''(x) = ax^ {\gamma},\, a>0,\, \gamma \in \mathbb R$, we give necessary and sufficient conditions for existence of $f$-minimal martingale measure.}

\section{Introduction}
\par This article is devoted to some important and exceptional properties of $f$-divergences. As known, the notion of $f$-divergence was introduced by Ciszar \cite{Ci} to measure the difference between two absolutely continuous probability measures by mean of the expectation of some convex function $f$ of their Radon-Nikodym density. More precisely, let $f$ be a convex function and $Z = \displaystyle\frac{dQ}{dP}$ be a Radon-Nikodym density of two measures $Q$ and $P$, $Q<\!< P$. Supposing that $f(Z)$ is integrable with respect to $P$,  $f$-divergence of $Q$ with respect to $P$ is defined as
$$f(Q|\!|P)= E_P[f(Z)].$$
One can remark immediately that this definition cover such important cases as variation distance when $f(x) = |x-1|$, as Hellinger distance when $f(x) = (\sqrt{x}-1)^ 2$ and Kulback-Leibler information when $f(x) = x \ln (x)$. Important as notion, $f$-divergence was studied in a number of books and articles (see for instance \cite{LV},  \cite{JSh})\\
\par In financial mathematics it is of particular interest to consider measures $Q^*$ which minimise on the set of all equivalent martingale measures the $f$-divergence. This fact is related to the introducing and studying so called incomplete models, like exponential Levy models, in which contingent claims cannot, in general, be replicated by admissible strategies. Therefore, it is important to determine strategies which are, in a certain sense optimal. Various criteria are used, some of which are linked to risk minimisation ( see \cite{FS}, \cite{S1}, \cite{S2}) and others consisting in maximizing certain utility functions (see \cite{GR}, \cite{K1}). It has been shown (see \cite{GR}, \cite{KSHA}) that such questions are strongly linked via Fenchel-Legendre transform to dual optimisation problems, namely to $f$-divergence minimisation on the set of equivalent martingale measures, i.e. the measures $Q$ which are equivalent to the initial physical measure $P$ and under which the stock price is a martingale. \\
\par Mentioned problems has been well studied in the case of relative entropy, when $f(x)=x\ln(x)$ (cf. \cite{M},\, \cite{FM}), also for power functions $f(x)=x^q$, $q>1$ or $q<0$ (cf. \cite{JKM}), $f(x)= - x^q$, $0<q<1$ (cf. \cite{CS}, \cite{CSL})  and for logarithmic divergence $f(x)=-\ln(x)$ (cf. \cite{Kl}), called common $f$-divergences.  Note that the three mentioned functions  all satisfy $f''(x)=ax^{\gamma}$ for an $a>0$ and a $\gamma\in\mathbb{R}$. The converse is also true, any function which satisfies $f''(x)=ax^{\gamma}$ is, up to linear term, a common $f$-divergence. It has in particular been noted that for these functions, the $f$-divergence minimal equivalent martingale measure, when it exists, preserves the Levy property, that is to say that the  law of Levy process under initial measure $P$ remains a law of Levy process under the $f$-divergence minimal equivalent martingale measure $Q^*$.\\
\par The aim of this paper is to study the questions of preservation of Levy property and associated properties such as  scaling property and invariance in time property for $f$-divergence minimal martingale measures when $P$ is a law of $d$-dimensional Levy process $X$ and $Q^*$ belongs to the set of so called equivalent martingale measures for exponential Levy model, i.e. measures under which the exponential of $X$ is a martingale. More precisely, let fix a convex function $f$ defined on $\mathbb R^{+,* }$ and  denote by $\mathcal M$ the set of equivalent martingale measures associated with exponential Levy model related to $X$. We recall that
an equivalent martingale measure $Q^*$ is $f$-divergence minimal if  $f(Z^*)$ is integrable with respect to $P$ where $Z^*$ is  the Radon-Nikodym density of $Q^*$ with respect to $P$, and
$$f(Q^*|\!|P)=\min_{Q\in \mathcal M }f(Q|\!|P).$$
We say that $Q^*$ preserves Levy property if $X$ remains Levy process under $Q^ *$.
The measure $Q^*$ is said to be scale invariant if  for all $x\in \mathbb R^+$,  $E_P|f(xZ^*)|<\infty$  and
$$f(xQ^*|\!|P)=\min_{Q\in \mathcal M }f(xQ|\!|P).$$
We also recall that an equivalent martingale measure $Q^*$ is said to be time invariant if for all $T>0$, and the restrictions $Q_T$,$P_T$ of the measures $P,Q$ on time interval $[0,T]$,  $E_P|f(Z_T^*)|<\infty$   and  
$$f(Q^*_T|\!|P_T)=\min_{Q\in \mathcal M }f(Q_T|\!|P_T)$$
\par In this paper  we study the shape of $f$  belonging  to the class of strictly convex tree times continuously differentiable functions and ones used as $f$-divergence, gives an equivalent martingale measure which preserves Levy property. More precisely, we consider  equivalent martingale measures $Q$ belonging to the class $\mathcal K ^ *$ such that for  all compact sets $K$ of $\mathbb R^{+,*}$
$$E_P|f(\frac{dQ_T}{dP_T})|<+\infty,\,\,\,\,E_Q|f'(\frac{dQ_T}{dP_T})|<+\infty,\,\,\,\,\sup_{t\leq T}\sup_{\lambda\in K}E_Q[f''(\lambda\frac{dQ_t}{dP_t})\frac{dQ_t}{dP_t}]<+\infty .$$

\par We denote by $Z_T^*$ Radon-Nikodym density of $Q_T^*$ with respect to $P_T$ and by $\beta ^*$ and $Y^ *$ the corresponding Girsanov parameters of an $f$-divergence minimal measure $Q^*$ on $[0,T]$, which preserves the Levy property and belongs to $\mathcal K^*$.

To precise the shape of $f$ we obtain  fundamental equations which necessarily verify $f$. Namely,  in the case $\stackrel{\circ}{supp}(\nu)\neq \emptyset$, for a.e. $x\in supp(Z^*_T)$ and a.e. $y\in supp(\nu)$, we prove that 
\begin{equation}
f'(xY^*(y))-f'(x)=\Phi(x)\sum_{i=1}^d \alpha_i (e^{y_i}-1)
\end{equation}
where  $\Phi$ is a continuously differentiable function on $\stackrel{\circ}{\mbox{supp}}(Z_T^*)$ and
$y = \,^\top (y_1,y _2,\cdots y _d)$, $\alpha = \,^\top (\alpha_1,\alpha _2,\cdots \alpha _d)$ are  vectors of $\mathbb{R}^d$. Furthermore, if $c\neq 0$, for a.e. $x\in \mbox{supp}(Z^*_T)$ and a.e. $y\in supp(\nu)$, we get that 
\begin{equation}
f'(xY^*(y))-f'(x)=xf''(x)\sum_{i=1}^d \beta_i^{*}(e^{y_i}-1)+\sum_{j=1}^d V_j(e^{y_j}-1)
\end{equation}
where $\beta^* = ^ \top\!(\beta^ *_1, \cdots, \beta ^*_d)$ is a first Girsanov parameter and $V= ^{\top}(V_1, \cdots, V_d)$ is a vector which belongs to the kernel of the matrix $c$, i.e. $cV=0$.
\par Mentioned above equations permit us to precise the form of $f$. Namely, we prove that
if  the set $\{\ln Y^*(y),\, y\in supp(\nu)\}$ is of non-empty interior and it contains zero, then
there exists $a>0$ and $\gamma\in\mathbb{R}$ such that for all $x\in supp(Z^*_t)$, 
\begin{equation}\label{res1}
f''(x)=ax^{\gamma}.
\end{equation}
Taking in account the known results we conclude that  in considered  case  the relation (\ref{res1}) is necessary and sufficient condition for $f$-divergence minimal martingale measure to preserve Levy property. In addition, as we will see, such $f$-divergence minimal measure will be  also scale and time invariant.\\
\par In the case when  ${}^{\top}\beta^* c\beta^*\neq 0$ and  support of $\nu$ is nowhere dense but when there exists at least one $y\in supp (\nu )$ such that $\ln (Y^*(y))\neq 0$, we prove that there exist $n\in \mathbb N$, the real constants $b_i, \tilde{b}_i, 1 \leq i\leq n,$ and $\gamma\in \mathbb R$, $a>0$ such that
$$f''(x)=  ax^{\gamma} + x^{\gamma}\sum_{i=1}^n b_i (\ln(x))^i + \frac{1}{x} \sum_{i=1}^n \tilde{b}_i (\ln(x))^{i-1}$$
The case when ${}^{\top}\beta^* c\beta^*= 0$ and $supp(\nu)$ is nowhere dense,  is not considered in this paper, and from what we know,  form an open question.\\

 We underline once more the exceptional properties of the class of functions such that:  $$f''(x)=ax^{\gamma}$$ and called common $f$-divergences. This class of functions is exceptional in a sense that they verify also scale and time invariance properties for all Levy processes. As well known,  $Q^*$ does not always exist. For some functions, in particular $f(x)=x\ln(x)$, or for some power functions, some necessary and sufficient conditions of existence of a minimal measure have been given (cf.\cite{HS},\cite{JKM}). We will give a unified version of these results for all functions which satisfy $f''(x)=ax^{\gamma}$, $a>0$, $\gamma\in\mathbb{R}$. We give also an example to show that the preservation of Levy property can have place not only for the functions verifying $f''(x)= ax^{\gamma} $.\\

The paper is organized in the following way:
in \ref{2}. we recall some known facts about exponential Levy models  and $f$-divergence minimal equivalent martingale measures. In \ref{3}. we give some known useful for us facts about $f$-divergence minimal martingale measures. In \ref{4}. we obtain  fundamental equations for Levy preservation property (Theorem \ref{feq} ). In \ref{5}. we give the result about the shape of $f$ having Levy preservation property for $f$-divergence minimal martingale measure (Theorem \ref{final}). In \ref{6}. we study the common f-divergences, i.e. with $f$ verifying $f''(x)=ax^{\gamma}$, $a>0$. Their properties are given in Theorem \ref{exist}.
\section{Some facts about Exponential Levy models}\label{2}
Let us describe our model in more details. We assume the financial market consists  of a bank account $B$ whose value at time $t$ is 
$$B_t=B_0e^{rt},$$
where $r\geq 0$ is the interest rate which we assume to be constant. We also assume that there are $d\geq 1$ risky assets whose prices are described by a $d$-dimensional stochastic process $S= (S_t)_{t\geq 0}$, 
$$S_t= \,^\top\!(S_0^{(1)}e^{X_t^{(1)}},\cdots,S_0^{(d)}e^{X_t^{(d)}})$$
where $X = (X_t)_{t\geq 0}$ is a $d$-dimensional Levy process, $X_t = \,^\top\!(X_t^{(1)},\cdots,X_t^{(d)})$ and $S_0=\,^\top \!(S_0^{(1)},\cdots,S_0^{(d)})$. We recall that Levy processes form the class of  processes with stationary and independent increment and that the characteristic function of the law of $X_t$ is given by the Levy-Khintchine formula : for all $t\geq 0$, for all $u\in\mathbb{R}$, 
$$E[e^{i<u,X_t>}]=e^{t\psi(u)}$$ where  
$$\psi(u)=i<u,b>-\frac{1}{2}\, ^{\top}\!ucu+
\int_{\mathbb{R}^d}[e^{i<u,y>}-1-i<u, h(y)> ]\nu(dy)$$
where $b\in\mathbb{R}^d$, $c$ is a positive $d\times d$ symmetric matrix, $h$ is a truncation function and $\nu$ is a Levy measure , i.e. positive measure on $\mathbb{R}^{d}\setminus \{0\}$ which satisfies 
$$\int_{\mathbb{R}^d}(1\wedge |y|^2) \nu(dy)<+\infty.$$
The triplet $(b,c,\nu)$ entirely determines the law of the Levy process $X$, and is called the characteristic triplet of $X$. From now on, we will assume that the interest rate $r=0$ as this will simplify calculations and the more general case can be obtained by replacing the drift $b$ by $b-r$. We also assume for simplicity that $S_0 =1$. 

\par We will denote by $\mathcal{M}$ the set of all locally equivalent martingale measures :
$$\mathcal{M}=\{Q\stackrel{\small{loc}}{\sim }P,\, S \text{ is a martingale under }Q \}.$$
We will assume that this set is non-empty, which is equivalent to assuming the existence of $Q\stackrel{\small{loc}}{\sim }P$ such that  the drift of $S$ under $Q$  is equal to zero.
We consider our model on finite time interval $[0,T]$, $T>0$, and for this reason the distinction between locally equivalent martingale measures and equivalent martingale measures does not need to be made. We recall that the density $Z$ of any equivalent to $P$ measure  can be written in the form $Z=\mathcal{E}(M)$ where $\mathcal{E}$ denotes the Doleans-Dade exponential and $M=(M_t)_{t\geq 0}$ is a local martingale. It follows from Girsanov theorem theorem that there exist predictable functions $\beta =\, ^{\top}\!(\beta^{(1)}, \cdots \beta^{(d)})$ and $Y$ verifying the integrability conditions : for $t\geq 0$ ($P$-a.s.)
$$ \int_0^{t} {^\top}\beta _s c \beta _s ds < \infty , $$ 
$$ \int_0^t \int_{\mathbb{R}^{d}}|\,h(y)\,(Y_s(y)-1)\,|\nu^{X, P}(ds, dy) < \infty , $$
and such that
\begin{equation}\label{mart}
M_t= \sum_{i=1}^d \int_0^t \beta^{(i)}_s dX^{c, (i)}_s+\int_0^t \int_{\mathbb{R}^{d}}(Y_s(y)-1)(\mu^{X}-\nu^{X, P})(ds, dy)
\end{equation}
where  $\mu^{X}$ is a  jumps measure  of the process $X$ and $ \nu^{X, P}$ is its compensator with respect to $P$ and the natural filtration $\mathbb F$, $\nu^{X, P}(ds, dy)= ds \,\nu(dy)$ ( for more details see \cite{JSh}). We will refer to $(\beta,Y)$ as the Girsanov parameters of the change of measure from $P$ into $Q$. It is known from Grigelionis result \cite{GREG} that a semi-martingale is a process with independent increments under $Q$ if and only if their semi-martingale characteristics are deterministic, i.e. the Girsanov parameters do not depend on $\omega$, i.e. $\beta$ depends only on time $t$ and $Y$ depends on time and jump size $(t,x)$. Since Levy process is homogeneous process, it implies that $X$ will remain a Levy process under $Q$ if and only if there exists $\beta\in\mathbb{R}$ and a positive measurable function $Y$ such that for all $t\leq T$ and all $\omega$, $\beta_t(\omega)=\beta$ and $Y_t(\omega,y)=Y(y)$. \\
 We recall that if Levy property is preserved, $S$ will be a martingale under $Q$ if and only if 
\begin{equation}\label{drift}
b+\frac{1}{2}diag(c)+c\beta+\int_{\mathbb{R}^{d}}[(e^y-1)Y(y)-h(y)]\nu(dy)=0
\end{equation}
where $e^y$ is a vector with components $e^{y_i}, 1\leq i\leq d,$ and $y= ^\top\!(y_1,\cdots,y_d)$.
This follows again from Girsanov theorem and reflects the fact that under $Q$ the drift of $S$ is equal to zero.

\section{Properties of $f$-divergence minimal martingale measures}\label{3}
Here we consider a fixed strictly convex continuously differentiable on $\mathbb R^{+,*}$ function $f$ and a time interval $[0,T]$. We recall in this section a few known and useful results about $f$-divergence minimisation on the set of equivalent martingale measures. Let  $(\Omega,\mathcal{F},\mathbb{F},P)$ be a probability filtered space with the natural filtration  $\mathbb F=(\mathcal F_t)_{t\geq t}$ satisfying usual conditions and let $\mathcal M$ be the set of equivalent martingale measures. We denote by $Q_t$ , $P_t$ the restrictions of the measures $Q$, $P$ on $\mathcal F_t$. We introduce Radon-Nikodym density process $Z= (Z_t)_{t\geq 0}$ related to $Q$, an equivalent martingale measure, where for $t\geq 0$
$$ Z_t= \frac{dQ_t}{dP_t}.$$ 
We denote by $Z^*$   Radon-Nikodym density process related with $f$-divergence minimal equivalent martingale measure $Q^*$.
\begin{definition}
An equivalent martingale measure $Q^*$ is said to be $f$-divergence minimal on the time interval $[0,T]$ if  $E_P|f(Z_T^*)|<\infty$  and 
$$E_P[f(Z_T^*)]=\min_{Q\in\mathcal{M}}E_P[f(Z_T)]$$
where $\mathcal M $ is a class of locally equivalent martingale measures.
\end{definition}
 
Then we introduce the subset of equivalent martingale measures
\begin{equation}\label{integcd0}
\mathcal{K}=\{Q\in\mathcal{M}\,|\,\, E_P|f(Z_T)|<+\infty \text{ and }E_Q[|f'(Z_T)|]<+\infty.\}
\end{equation}
We will concentrate ourselves on the case when the minimal measure, if it exists, belongs to $\mathcal{K}$. Note that for a certain number of functions this is necessarily the case.
\begin{lemma}[cf \cite{LV}, Lemma 8.7]
Let $f$ be a convex continuously differentiable on $\mathbb R ^{+,* }$ function. Assume that for $c>1$ there exist positive constants $c_0,c_1,c_2,c_3$ such that for $u>c_0$,
\begin{equation}\label{1l}
f(cu)\leq c_1f(u)+c_2u+c_3
\end{equation}
Then a measure $Q\in\mathcal{M}$ which is $f$-divergence minimal necessarily belongs to $\mathcal{K}$. 
\end{lemma}
We now recall the following necessary and sufficient condition for a martingale measure to be minimal. 
\begin{theorem}[cf \cite{GR}, Th 2.2]\label{GR1}
Consider $Q^*\in\mathcal{K}$. Then,  $Q^*$ is minimal if and only if for all $Q\in\mathcal{K}$, 
$$E_{Q^*}[f'(Z^*_T)]\leq E_{Q}[f'(Z^*_T)].$$
\end{theorem}
This result is in fact true in the much wider context of semi-martingale modelling. We will mainly use it here to check that a candidate is indeed a minimal measure. We will also use extensively another result from \cite{GR} in order to obtain conditions that must be satisfied by minimal measures. 
\begin{theorem}[cf \cite{GR}, Th 3.1]\label{GR2}
Assume $Q^*\in \mathcal{K}$ is an $f$-divergence minimal martingale measure. Then there exists $x_0\in\mathbb{R}$ and a predictable $d$-dimensional process $\phi$ such that 
$$f'(\frac{dQ^*_T}{dP_T})=x_0 +\sum_{i=1}^d\int _0^T\phi _t^{(i)}\, d S_t^{(i)}$$
and such that $\sum_{i=1}^d\int _0^{\cdot}\phi _t^{(i)}\, d S_t^{(i)}$ defines a martingale under the measure $Q^*$.
\end{theorem}

\section{A fundamental equation for $f$-divergence minimal Levy preserving martingale measures}\label{4}
Our main aim in this section is to obtain an equation satisfied by the Radon-Nikodym density of $f$-divergence minimal equivalent martingale measures. This result will both enable us to obtain information about the Girsanov parameters of $f$-divergence minimal equivalent martingale measures and also to determine conditions which must be satisfied by the function $f$ in order to a f-minimal equivalent martingale measure exists. Let us introduce the class $\mathcal K ^ *$ of locally equivalent martingale measures verifying: for all compact sets $K$ of $\mathbb R^{+,*}$
\begin{equation}\label{integcd}
E_P|f(Z_T)|<+\infty,\,\,\,\,E_Q|f'(Z_T)|<+\infty,\,\,\,\,\sup_{t\leq T}\sup_{\lambda \in K}E_Q[f''(\lambda Z^*_t)Z^*_t]<+\infty.
\end{equation}

\begin{theorem}\label{feq}
Let $f$ be strictly convexe ${\mathcal C}^3(\mathbb R^{+,*})$ function. Let $Z^*$ be the density of an $f$-divergence minimal measure $Q^*$ on $[0,T]$, which preserves the Levy property and belongs to $\mathcal K^*$.
We denote by $(\beta^*,Y^*)$ its Girsanov parameters. Then, if $\stackrel{\circ}{supp}(\nu)\neq \emptyset \, $, for a.e. $x\in supp(Z^*_T)$ and a.e. $y\in supp(\nu)$, we have 
\begin{equation}\label{eqfmin}
f'(xY^*(y))-f'(x)=\Phi(x)\sum_{i=1}^d \alpha_i (e^{y_i}-1)
\end{equation}
where $\Phi$ is a continuously differentiable function on $\stackrel{\circ}{\mbox{supp}}(Z_T^*)$ and $\alpha = \,^\top (\alpha_1,\alpha _2,\cdots \alpha _d)$ is a vector of $\mathbb{R}^d$. Furthermore, if $c\neq 0$, for a.e. $x\in \mbox{supp}(Z^*_T)$ and a.e. $y\in supp(\nu)$, we have 
\begin{equation}\label{equ}
f'(xY^*(y))-f'(x)=xf''(x)\sum_{i=1}^d \beta_i^{*}(e^{y_i}-1)-\sum_{j=1}^d V_j(e^{y_j}-1)
\end{equation}
where $\beta^* = ^ \top\!(\beta^ *_1, \cdots, \beta ^*_d)$ and $V= ^{\top}(V_1, \cdots, V_d)$ belongs to the kernel of the matrix $c$, i.e. $cV=0$.
\end{theorem}

We recall that for all $t\leq T$, since $Q^*$ preserves Levy property, $Z^*_t$ and $\displaystyle\frac{Z^*_T}{Z^*_t}$ are independent under $P$ and that $\mathcal{L}(\displaystyle\frac{Z^*_T}{Z^*_t})=\mathcal{L}(Z^*_{T-t})$.  Therefore denoting
 $$\rho(t,x)=E_{Q^*}[f'(x Z^*_{T-t})],$$
and taking cadlag versions of processes, we deduce that $Q^*$-a.s. for all $t\leq T$
$$E_{Q^*}[ f'( Z^*_T)|\mathcal{F}_t]=\rho(t,Z^*_t)$$

We note that the proof of Theorem \ref{feq} is based on the identification using Theorem \ref{GR2} and an application of decomposition formula to function $\rho$. However, the function $\rho$
is not necessarily twice continuously differentiable in $x$ and once continuously differentiable in $t$. So, we will proceed by approximations, by application of Ito formula to specially constructed function $\rho _n$. In order to do this, we need  a number of auxiliary lemmas given in the next section.

Since the result of Theorem \ref{feq} is strongly related to the support of $Z^*_T$, we are also interested with the question : when this support is  an interval? This question has been well studied in \cite{T}, \cite{Sa} for infinitely divisible distributions. In our case, the specific form of the Girsanov parameters following from preservation of Levy property allow us to obtain the following result proved in subsection 4.3. 

\begin{proposition}\label{support}
Let $Z^*$ be the density of an $f$-divergence minimal equivalent martingale measure on $[0,T]$, which preserves the Levy property and belongs to $\mathcal K ^*$. Then 
\begin{enumerate}
\item[(i)] If ${}^{\top}\beta^* c\beta^*\neq 0$, then $supp(Z^*_T)=\mathbb{R}^{+,*}$. 
\item[(ii)] If ${}^{\top}\beta^* c\beta^*= 0$, $\stackrel{\circ}{supp}(\nu)\neq \emptyset$, $0\in supp(\nu)$ and $Y^ *$ is not identically 1 on $\stackrel{\circ}{supp}~(\nu)$ , then
\begin{enumerate}
\item[(j)] in the case $\ln (Y(y))>0$ for all $y\in supp (\nu)$, there exists $A>0$ such that  $ supp(Z^*_T)=[A,+\infty[$; 
\item[(jj)] in the case  $\ln (Y(y))<0$ for all $y\in supp (\nu)$ there exists $A>0$ such that  $supp(Z^*_T)=]0,A]$;
\item[(jjj)] in the case when there exist $y, \bar{y}\in supp(\nu)$ such that $\ln(Y^*(y)).\ln(Y^*(\bar{y}))<0$, we have $supp(Z^*_T)=\mathbb{R}^{+,*}$.
\end{enumerate} 
\end{enumerate}
\end{proposition}

\subsection{Some auxiliary lemmas}
We begin with approximation lemma. Let a strictly convex tree times continuously differentiable on $ \mathbb R ^{+,*}$ function $f$ be fixed.

\begin{lemma}\label{approx} There exists a sequence of bounded functions $(\phi_n)_{n\geq 1}$, which are of class $\mathcal{C}^2$ on $\mathbb{R}^{+*}$, increasing, such that for all $n\geq 1$, $\phi_n$ coincides with $f'$ on the compact set $[\frac{1}{n},n]$ and such that for sufficiently big $n$ the following inequalities hold for all $x,y>0$~:
\begin{equation}\label{apprfprime}
|\phi_n(x)|\leq 4|f'(x)|+ \alpha \text{ , }|\phi'_n(x)|\leq 3f''(x) \text{ , }|\phi_n(x)-\phi_n(y)|\leq 5|f'(x)-f'(y)|
\end{equation}
where $\alpha$ is a real positive constant.
\end{lemma} 

\noindent \it Proof \rm 
We set, for $n\geq 1$, 
$$A_n(x)=f'(\frac{1}{n})-\int_{x\vee \frac{1}{2n} }^{\frac{1}{n}}f''(y)(2ny-1)^2(5-4ny)dy$$
$$B_n(x)=f'(n)+\int_n^{x\wedge (n+1)} f''(y)(n+1-y)^2(1+2y-2n)dy$$
and finally
$$\phi_n(x)=\begin{cases}
& A_n(x) \text{ if }0 \leq x< \frac{1}{n},
\\&f'(x) \text{ if }\frac{1}{n}\leq x\leq n,
\\& B_n(x) \text{ if }x>n.

\end{cases}$$

Here $A_n$ and $B_n$ are defined so that $\phi_n$ is of class $\mathcal{C}^2$ on $\mathbb{R}^{+,*}$.
For the inequalities we use the fact that $f'$ is increasing function as well as the estimations:
$0\leq (2nx-1)^2(5-4nx)\leq 1$ for $\frac{1}{2n} \leq x\leq \frac{1}{n}$ and
$0\leq(n+1-x)^2(1+2x-2n)\leq 3$ for $n \leq x\leq n+1$. 
$\Box$

\vskip 0.2cm Let $Q$ be Levy property preserving locally equivalent martingale measure and $(\beta, Y)$ its Girsanov parameters when change from $P$ into $Q$. We use  the function 
$$\rho_n(t,x)=E_Q[\phi_n( xZ_{T-t})]$$
to obtain the following analog to Theorem \ref{thdec}, replacing  $f'$ with $\phi_n$. 

For this let us denote for $0\leq t \leq T$
\begin{equation}\label{xin}
\xi_t^{(n)}(x)= E_Q[\phi'_n(xZ_{T-t})\,Z_{T-t}]
\end{equation}
and
\begin{equation}\label{han}
H_t^{(n)}(x,y)=E_Q[\phi_n(xZ_{T-t}Y(y))-\phi_n(xZ_{T-t})]
\end{equation}

\begin{lemma}\label{decomp}
We have $Q^*$-a.s., for all $t\leq T$,
\begin{equation}\label{ndecomp}
\rho_n(t,Z_t)=E_Q[\phi_n( Z_T)] +
\end{equation}
$$\sum_{i=1}^d\beta _i \int_0^t  \xi^{(n)}_s ( Z_{s-})\,Z_{s-} dX^{(c),Q,i}_s+\int_0^t \int_{\mathbb{R}^{d}}H^{(n)}_s( Z_{s-},y)\,(\mu^X-\nu ^{X,Q})(ds, dy)$$
where $\beta = \,^\top\!(\beta _1,\cdots, \beta _d) $ and $\nu ^{X,Q}$ is a compensator of the jump measure $\mu^X$ with respect to $(\mathbb F, Q)$.
\end{lemma}
\noindent \it Proof \rm In order to apply the Ito formula to $\rho_n$, we need to show that $\rho_n$ is twice continuously differentiable with respect to $x$ and once with respect to $t$ and that the corresponding derivatives are bounded for all $t\in[0,T]$ and $x\geq \epsilon$, $\epsilon >0.$ First of all, we note that from the definition of $\phi_n$
for all $x\geq \epsilon >0$
$$|\frac{\partial}{\partial x}\phi_n( xZ_{T-t})| = | Z_{T-t}\phi'_n( xZ_{T-t})|\leq \frac{(n+1)}{\epsilon}\sup_{z>0}|\phi'_n(z)|<+\infty.$$
Therefore, $\rho_n$ is differentiable with respect to $x$ and we have
$$\frac{\partial}{\partial x}\rho_n(t,x)= E_Q[\phi'_n( xZ_{T-t})\,Z_{T-t}].$$
Moreover, the function $(x,t)\mapsto \phi'_n( xZ_{T-t})Z_{T-t}$ is continuous $P$-a.s. and bounded. This implies that $\frac{\partial}{\partial x}\rho_n$ is continuous and bounded for
$t\in [0,T] $ and $x\geq \epsilon $. \\
 In the same way, for all $x\geq \epsilon >0$
$$|\frac{\partial^2}{\partial x^2}\phi_n( xZ_{T-t})|= Z_{T-t}^2\phi''_n( x Z_{T-t})\leq \frac{(n+1)^2}{\epsilon ^ 2} \sup_{z>0}\phi''_n(z)<+\infty .$$
Therefore, $\rho_n$ is twice continuously differentiable in $x$ and 
$$\frac{\partial^2}{\partial x^2}\rho_n(t,x)= E_Q[\phi_n''( xZ_{T-t})Z^2_{T-t}]$$
We can verify easily that it is again continuous and bounded function.
In order to obtain differentiability with respect to $t$, we need to apply the Ito formula to $\phi_n$ :
$$\begin{aligned}
\phi_n( xZ_{t})=&\phi_n( x)+\sum_{i=1}^d \int_0^{t} x\phi'_n( xZ_{s-})\beta _iZ_{s-}dX^{(c),Q,i}_s
\\&+\int_0^{t}\int_{\mathbb{R}^{d}}\phi_n( xZ_{s-}Y(y))-\phi_n( xZ_{s-})\,(\mu^X-\nu^{X,Q})(ds, dy)
\\&+\int_0^{t}\psi_n( x,Z_{s-})ds
\end{aligned}$$
where 
$$\begin{aligned}
\psi_n( x,Z_{s-})=& ^{\top}\beta c\beta [ x Z_{s-}\phi'_n( xZ_{s-})+\frac{1}{2}\,x^2 \,Z_{s-}^ 2 \phi''_n( xZ_{s-})]
\\&+\int_{\mathbb{R}^{d}}[(\phi_n(  xZ_{s-}Y(y))-\phi_n( xZ_{s-}))\,Y(y)- x\phi'_n( xZ_{s-})Z_{s-}(Y(y)-1)]\nu(dy).
\end{aligned}$$
Therefore, 
$$E_Q[\phi_n( xZ_{T-t})]=\int_0^{T-t}E_Q[\psi_n( x,Z_{s-})]ds$$
so that $\rho_n$ is differentiable with respect to $t$ and 
$$\frac{\partial}{\partial t}\rho_n(t,x)=-E_Q[\psi_n( x,Z_{s-})]_{|_{s=(T-t)}}$$
We can also easily verify that this function is continuous and bounded. For this we take in account the fact that $ \phi_n$, $\phi'_n$ and  $\phi''_n $  are bounded functions and also that the Hellinger process
of $Q_T$ and $P_T $ of the order $1/2$  is finite.
\\ 

We can finally apply the Ito formula to $\rho_n$. For that we use the stopping times
$$s_m= \inf \{ t\geq 0\,|\, Z_t \leq \frac{1}{m} \},$$
with $m\geq 1$ and $\inf\{\emptyset\}= +\infty$.
Then, from Markov property of L\'evy process we have~:
$$\rho _n(t\wedge s_m, Z_{t\wedge s_m}) = E_Q(\phi _n(\lambda Z_T)\,|\, \mathcal F_{t\wedge s_m})$$
We remark that $(E_Q(\phi _n(\lambda Z_T)\,|\, \mathcal F_{t\wedge s_m})_{t\geq 0}$ is $Q$-martingale, uniformly integrable with respect to $m$. From Ito formula we have :
\begin{eqnarray*} \rho _n(t\wedge s_m, Z_{t\wedge s_m}) = E_Q(\phi _n(\lambda Z_T)) + \int_0^{t\wedge s_m}\frac{\partial \rho _n}{\partial s}(s, Z_{s-}) ds\\\\+ \int_0^{t\wedge s_m}\frac{\partial \rho _n}{\partial x}(s, Z_{s-}) dZ_s +
 \frac{1}{2}\int_0^{t\wedge s_m}\frac{\partial^2 \rho _n}{\partial x^2}(s, Z_{s-}) d< Z^c>_s \\\\ +
\sum _{0\leq s\leq t\wedge s_m}\rho _n(s, Z_{s})- \rho _n(s, Z_{s-}) - \frac{\partial \rho _n}{\partial x}(s, Z_{s-}) \Delta Z_s
\end{eqnarray*}
where $\Delta Z_s = Z_s - Z_{s-}$.
After some standard simplifications, we see that
$$\rho _n(t\wedge s_m, Z_{t\wedge s_m}) =   A_{t\wedge s_m} + M_{t\wedge s_m}$$
where $(A_{t\wedge s_m})_{0\leq t\leq T}$ is predictable process, which is equal to zero,
\begin{eqnarray*}A_{t\wedge s_m} = \int_0^{t\wedge s_m}\frac{\partial \rho _n}{\partial s}(s, Z_{s-}) ds +
 \frac{1}{2}\int_0^{t\wedge s_m}\frac{\partial^2 \rho _n}{\partial x^2}(s, Z_{s-}) d< Z^c>_s +\\\\
\int_0^{t\wedge s_m}\int _{\mathbb R}[ \rho _n(s, Z_{s-}+x)- \rho _n(s, Z_{s-}) - \frac{\partial \rho _n}{\partial x}(s, Z_{s-})x] \nu^ {Z,Q}(ds, dx)
\end{eqnarray*}
and $(M_{t\wedge s_m})_{0\leq t\leq T}$ is a $Q$-martingale,
\begin{eqnarray*}M_{t\wedge s_m} =  E_Q(\phi _n(\lambda Z_T)) +  \int_0^{t\wedge s_m}\frac{\partial \rho _n}{\partial x}(s, Z_{s-}) dZ^c_s+\\
\int_0^{t\wedge s_m}\int _{\mathbb R}[ \rho _n(s, Z_{s-}+x)- \rho _n(s, Z_{s-})](\mu^ Z(ds, dx)-\nu^ {Z,Q}(ds, dx))
\end{eqnarray*}
Then, we pass to the limit as $m\rightarrow +\infty$. We remark that the sequence $(s_m)_{m\geq 1}$ is going to $+\infty$ as $m\rightarrow \infty$. From \cite{RY}, corollary 2.4, p.59, we obtain that
$$\lim_{m\rightarrow \infty} E_Q(\phi_n(Z_T)\,|\, \mathcal F_{t\wedge s_m})= E_Q(\phi_n(Z_T)\,|\, \mathcal F_{t})$$ and by the definition of local martingales we get:
$$\lim_{m\rightarrow \infty} \int_0^{t\wedge s_m}\frac{\partial \rho _n}{\partial x}(s, Z_{s-}) dZ^c_s =  \int_0^{t}\frac{\partial \rho _n}{\partial x}(s, Z_{s-}) dZ^c_s = \int_0^{t}\lambda \xi_s^{(n)}(Z_{s-}) dZ^c_s$$
and
$$\lim_{m\rightarrow \infty}\int_0^{t\wedge s_m}\int _{\mathbb R}[ \rho _n(s, Z_{s-}+x)- \rho _n(s, Z_{s-})](\mu^ Z(ds, dx)-\nu^ {Z,Q}(ds, dx)) =$$
$$ \int_0^{t}\int _{\mathbb R}[ \rho _n(s, Z_{s-}+x)-\rho _n(s, Z_{s-})](\mu^ Z(ds, dx)-\nu^ {Z,Q}(ds, dx))$$
Now, in each stochastic integral we pass from the integration with respect to the process $Z$ to the one with respect to the process $X$. For that we remark that
$$dZ^c_s = \sum _{i=1}^d \beta ^{(i)} Z_{s-} dX_s^{c,Q,i},\,\,\,\Delta Z_s = Z_{s-}Y(\Delta X_s).$$ Lemma \ref{decomp} is proved.
$\Box $

\subsection{A decomposition for the density of Levy preserving martingale measures}
This decomposition will follow from a previous one by a limit passage. Let again $Q$ be Levy property preserving locally equivalent martingale measure and $(\beta, Y)$ the corresponding Girsanov parameters when passing from $P$ to $Q$. We introduce
cadlag versions of the following processes: for $t>0$
$$\xi_t(x)=  E_Q[f''(xZ_{T-t})Z_{T-t}]$$
and
\begin{equation}\label{H}
 H_t(x,y)=E_Q[f'(xZ_{T-t}Y(y))-f'(xZ_{T-t})]
\end{equation}
\begin{theorem}\label{thdec}
Let $Z$ be the density of a Levy preserving equivalent martingale measure $Q$. Assume that 
$Q$ belongs to $\mathcal K ^*$.
 Then we have $Q$- a.s, for all $t\leq T$,  
\begin{equation}\label{decf}
E_Q[f'( Z_T)|\mathcal{F}_t]=E_Q[f'( Z_T)] +
\end{equation}
$$\sum_{i=1}^d \beta _i\int_0^t \xi_s( Z_{s-}) Z_{s-} d {X}^{(c),Q,i}_s+\int_0^t \int_{\mathbb{R}^d}H_s( Z_{s-},y)\,(\mu^{X}-\nu^{X, Q})(ds, dy)$$

\end{theorem}

We now turn to the proof of Theorem \ref{thdec}. In order to obtain the decomposition for $f'$, we obtain convergence in probability of the different stochastic integrals appearing in Lemma \ref{decomp}.
\par \it Proof of Theorem \ref{thdec} \rm For a  $n\geq 1$, we introduce the stopping times 
\begin{equation}\label{taun}
\tau_n= \inf\{t\geq 0\,|\, Z_t\geq n \text{ or } Z_t\leq \frac{1}{n}\}
\end{equation}
where $\inf\{\emptyset\}= +\infty$ and we note that $\tau _n\rightarrow +\infty$ ($P$-a.s.) as $n\rightarrow \infty\,$.
First of all, we note that 
$$|E_Q[f'( Z_T)|\mathcal{F}_t]-\rho_n(t,Z_t)|\leq E_Q[|f'( Z_T)-\phi_n( Z_T)||\mathcal{F}_t]$$
As $f'$ and $\phi_n$ coincide on the interval $[\frac{1}{n},n]$, it follows from Lemma \ref{decomp} that 
$$\begin{aligned}
|E_Q[f'( Z_T)|\mathcal{F}_t]-\rho_n(t,Z_t)|&\leq E_Q[|f'( Z_T)-\phi_n( Z_T)|{\bf 1}_{\{\tau_n\leq T\}}|\mathcal{F}_t]
\\&\leq  E_Q[(5|f'( Z_T)|+ \alpha){\bf 1}_{\{\tau_n\leq T\}}|\mathcal{F}_t].
\end{aligned}$$
Now, for every $\epsilon >0$, by Doob inequality and Lebesgue dominated convergence theorem we get: 
$$\lim_{n\to+\infty}Q(\sup_{t\leq T}E_Q[(5|f'( Z_T)|+\alpha){\bf 1}_{\{\tau_n\leq T\}}|\mathcal{F}_t] > \epsilon )\leq \lim_{n\to+\infty}\frac{1}{\epsilon }E_Q[(5|f'( Z_T)|+\alpha){\bf 1}_{\{\tau_n\leq T\}}]=0$$
Therefore, we have 
$$\lim_{n\to+\infty}Q(\sup_{t\leq T}|E_Q[f'( Z_T) - \rho_n(t, Z_t)|\mathcal{F}_t]|>\epsilon)=0.$$
We now turn to the convergence of the three elements of the right-hand side of (\ref{ndecomp}). 
We have almost surely $\lim_{n\to+\infty}\phi_n( Z_T)=f'( Z_T)$, and for all $n\geq 1$, $|\phi_n( Z_T)|\leq 4|f'( Z_T)|+\alpha $. Therefore, it follows from the dominated convergence theorem that, 
$$\lim_{n\to+\infty}E_Q[\phi_n( Z_T)]=E_Q[f'( Z_T)].$$
 We prove now  the convergence of continuous martingale parts of (\ref{ndecomp}). It follows from Lemma \ref{approx} that 
$$\begin{aligned}                       
Z_t\,|\xi^{(n)}_t( Z_t)-\xi_t( Z_t)|\leq &  E_Q[Z_T|\phi'_n( Z_T)-f''( Z_T)|\,|\,\mathcal{F}_t]\leq
\\& 4 E_Q[Z_T|f''( Z_T)|{\bf 1}_{\{\tau_n\leq T\}}|\mathcal{F}_t].
\end{aligned}$$                     
Hence, we have as before  for $\epsilon >0 $           
$$\lim_{n\to+\infty}Q(\sup_{t\leq T}Z_t\,|\xi^{(n)}_t( Z_t)-\xi_t( Z_t)|>\epsilon)\leq \lim_{n\to+\infty}\frac{4}{\epsilon}E_Q[Z_Tf''( Z_T){\bf 1}_{\{\tau_n\leq T\}}]=0$$

 Therefore, it follows from the Lebesgue dominated convergence theorem for stochastic integrals (see \cite{JSh}, Theorem I.4.31, p.46 ) that for all $\epsilon >0$ and $1\leq i\leq d$ 
$$\lim_{n\to+\infty}Q(\sup_{t\leq T}\,\big|\int_0^t Z_{s-}\,(\xi^{(n)}_s( Z_{s-})-\xi_s( Z_{s-}))dX^{(c),Q,i}_s \big|>\epsilon)=0.$$
It remains to show the convergence of the discontinuous martingales to zero as $n\rightarrow~\infty $.
 We start by writing 
$$\int_0^{t}\int_{\mathbb{R}^{d}}[H^{(n)}_s( Z_{s-},y)-H_s( Z_{s-},y)](\mu^X-\nu^{X, Q})(ds, dy) =M^{(n)}_t+N^{(n)}_t$$
with 
$$M^{(n)}_t=\int_0^t \int_{ \mathcal A}[H^{(n)}_s( Z_{s-},y)-H_s( Z_{s-},y)] (\mu^X-\nu^{X, Q})(ds, dy),$$
$$N^{(n)}_t=\int_0^t \int_{\mathcal A^ c}[H^{(n)}_s( Z_{s-},y)-H_s( Z_{s-},y)] (\mu^X-\nu^{X, Q})(ds, dy),$$
where $\mathcal A = \{y : |Y(y)-1|<\frac{1}{4}\}$.\\
For $p\geq 1$, we consider the sequence of stopping times $\tau_p$ defined by (\ref{taun})   with replacing $n$ by real positive $p$. We introduce also the  processes $$M^{(n,p)}= (M^{(n,p)}_t)_{t\geq 0},\,\,
N^{(n,p)}= (N^{(n,p)}_t)_{t\geq 0}$$ with $M^{(n,p)}_t= M^{(n)}_{t\wedge \tau _p}$, $N^{(n,p)}_t= N^{(n)}_{t\wedge \tau _p}$.
We remark that for $p\geq 1$ and $\epsilon >0$
$$ Q(\sup _{t\leq T} |M^{(n)}_t+N^{(n)}_t| >\epsilon) \leq Q(\tau _p \leq T) + Q( \sup _{t\leq T} |M^{(n,p)}_t| >\frac{\epsilon }{2}) +
 Q( \sup _{t\leq T} |N^{(n,p)}_t| >\frac{\epsilon }{2}).$$
Furthermore, we obtain  from Doob martingale inequalities that
\begin{equation}\label{02}
 Q( \sup _{t\leq T} |M^{(n,p)}_t|  >\frac{\epsilon}{2}) \leq \frac{4}{\epsilon ^ 2}\mathbb E_Q [(M^{(n,p)}_T)^ 2]
\end{equation}
and
\begin{equation}\label{03}
 Q( \sup _{t\leq T} |N^{(n,p)}_t|  >\frac{\epsilon }{2}) \leq \frac{2}{\epsilon }\mathbb E_Q |N^{(n,p)}_T|
\end{equation} 
Since $\tau _p \rightarrow + \infty$ as $p\rightarrow +\infty$ it is sufficient to show  that $E_Q[M_T^{(n,p)}]^2$ and $E_Q|N_T^{(n,p)}|$ converge to $0$ as $n\rightarrow \infty $.\\
For that we estimate $E_Q [(M^{(n,p)}_T)^ 2]$ and prove that\\\\ 
$E_Q [(M_T^{(n,p)})^ 2]\leq$\\
$$C \big(\int _0^ T\sup _{v\in K}\mathbb E_Q^2[ Z_{s}\,f''( v Z_{s}){\bf 1}_{\{\tau _{q_n}<s\}}]ds\big)\,\big(\int _ {\mathcal A} (\sqrt{Y(y)}-1)^ 2\nu(dy) \big)$$
where $C$ is a constant, $K$ is some compact set of $\mathbb R^{+,*}$ and $q_n= \frac{ n }{4p }$.\\\\
First we note  that on stochastic interval $[\![0, T\wedge \tau _p)]\!]$ we have $1/p\leq  Z_{s-}\leq p$, and, hence,
$$E_Q[(M_T^{(n,p)})^2]=E_Q[\int_0^{T\wedge \tau _p}\int_{\mathcal A}|H^{(n)}_s( Z_{s-},y)-H_s( Z_{s-},y)|^2\,Y(y)\nu(dy)ds ]\leq$$
$$\int_0^T\int_{\mathcal A}\sup_{1/p\leq x\leq p}|H^{(n)}_{T-s}( x,y)-H_{T-s}( x,y)|^2\,Y(y)\nu(dy)ds $$
To estimate the difference $|H^{(n)}_{T-s}( x,y)-H_{T-s}( x,y)|$ we note that
$$H^{(n)}_{T-s}( x,y)-H_{T-s}( x,y)= E_Q[\phi_n( xZ_{s}Y(y))-\phi_n( xZ_{s})-f'( xZ_{s}Y(y))+f'( xZ_{s})]$$
From  Lemma \ref{approx} we deduce that if $ x Z_{s} Y(y) \in [1/n, n]$ and $ x Z_{s} \in [1/n, n]$ then the expression on the right-hand side of the previous expression is zero. But if $y\in \mathcal A$ we also have : $3/4 \leq Y(y)\leq 5/4$ and, hence,\\
$|H^{(n)}_{T-s}( x,y)-H_{T-s}( x,y)|\leq$ \\ 
$$|E_Q[{\bf 1}_{\{\tau _{q_n}\leq s\}}|\phi_n( xZ_{s}Y(y))-\phi_n( xZ_{s})-f'( xZ_{s}Y(y))+f'( xZ_{s})|]$$
Again from the inequalities of Lemma \ref{approx} we get:
$$|H^{(n)}_{T-s}( x,y)-H_{T-s}( x,y)|\leq  6 E_Q[{\bf 1}_{\{\tau _{q_n}\leq s\}}|f'( xZ_{s}Y(y))-f'( xZ_{s})|]$$
Writing 
$$f'( xZ_{s}Y(y))-f'( xZ_{s})= \int _1^ {Y(y)} x Z_{s} f''( x Z_{s}\theta) d\theta$$
we finally get
$$|H^{(n)}_{T-s}( x,y)-H_{T-s}( x,y)|\leq  6\sup_{3/4\leq u \leq 5/4}E_Q[{\bf 1}_{\{\tau _{q_n}\leq s\}}\, x Z_{s}\,f''( x uZ_{s})]|Y(y)-1|$$
and this gives us the estimation of $E_Q [(M^{(n,p)}_T)^ 2]$ cited above.\\
We know that $P_T\sim Q_T$ and this means that the corresponding Hellinger process of order 1/2 is finite:
$$ h_T(P, Q, \frac{1}{2}) =\frac{T}{2}\,^{\top}\beta c \beta  + \frac{T}{8} \int_{\mathbb R}(\sqrt{Y(y)}-1)^ 2\nu(dy) < +\infty .$$
Then 
$$\int_{\mathcal A}(\sqrt{Y(y)}-1)^ 2\nu(dy) < +\infty .$$
From Lebesgue dominated convergence theorem  and (\ref{integcd}) we get:
$$\int _0^ T\sup _{v\in K}\mathbb E_Q^2[ Z_{s}f''( v Z_{s}){\bf 1}_{\{\tau _{q_n}\leq s\}}]ds \rightarrow 0$$
as $n\rightarrow +\infty$ and this information together with the estimation of $E_Q[(M_T^{(n,p)})^2]$ proves the convergence of $E_Q[(M_T^{(n,p)})^2]$  to zero as $n\rightarrow +\infty$.

\par We now turn to the convergence of $E_Q| N_T^{(n,p)}|$ to zero as $n\rightarrow +\infty$. 
For this we prove that
$$E_Q| N_T^{(n,p)}|\leq 2 T E_Q[{\bf 1}_{\{\tau _{n}\leq T\}}(5|f'( Z_{T})|+\alpha)]\int_{\mathcal A^ c}Y(y) d\nu $$
We start by noticing that
$$E_Q|N^{(n,p)}_T |\leq 2 E_Q[\int_0^{T\wedge \tau _p}\int_{\mathcal A^c}|H^{(n)}_s( Z_{s-},y)-H_s( Z_{s-},y)|\,Y(y)\nu(dy)ds]\leq $$
$$ 2 \int_0^{T}\int_{\mathcal A^c}E_Q[|H^{(n)}_s( Z_{s-},y)-H_s( Z_{s-},y)|\,Y(y)\nu(dy)ds]$$
 To evaluate the right-hand side of previous inequality we write\\
$|H^{(n)}_s( x,y)-H_s( x,y)|\leq$\\
 $$E_Q|\phi_n( xZ_{T-s}Y(y))-f'( xZ_{T-s}Y(y))|  + E_Q |\phi_n( xZ_{T-s})-f'( xZ_{T-s})|.$$
We remark that in law with respect to $Q$
$$ |\phi_n( x Z_{T-s}Y(y))-f'( x Z_{T-s}Y(y))| = E_Q[|\phi_n( Z_{T})-f'( Z_{T})|\,|\,Z_{s}=x\,Y(y)]$$
and
$$ |\phi_n( xZ_{T-s})-f'( xZ_{T-s})| = E_Q[|\phi_n( Z_{T})-f'( Z_{T})|\,|\,Z_{s}=x]$$
Then
$$H^{(n)}_s( x,y)-H_s( x,y)| \leq 2 E_Q|\phi_n( Z_{T})-f'( Z_{T})|$$
From  Lemma \ref{approx} we get:
$$E_Q|\phi_n( xZ_{T})-f'( xZ_{T})|\leq E_Q[{\bf 1}_{\{\tau _{n}\leq T\}}|\phi_n( Z_{T})-f'( Z_{T})|]\leq E_Q[{\bf 1}_{\{\tau _{n}\leq T\}}(5|f'( Z_{T})|+\alpha)]$$
and is proves the estimation for $E_Q|N^{(n,p)}_T |$.

Then, Lebesgue dominated convergence theorem applied for the right-hand side of the previous inequality shows that it tends to zero as $n\rightarrow \infty$.
On the other hand, from the fact that the Hellinger process is finite and also from the inequality $(\sqrt{Y(y)}-1)^2\geq Y(y)/25$ verifying on $\mathcal A^ c$
we get
$$\int_{A^ c}Y(y) d\nu < +\infty$$
This result with previous convergence  prove the convergence of $E_Q|N^{(n,p)}_T |$ to zero as $n\rightarrow \infty$. Theorem 4 is proved.
$\Box$. 

\subsection{Proof of Theorem \ref{feq} and Proposition \ref{support}}
\noindent \it Proof of Theorem \ref{feq}. \rm We define a process $\hat{X}= \,^\top(\hat{X}^{(1)}, \cdots \hat{X}^{(d)})$ such that for $1\leq i\leq d$ 
and $t\in[0,T]$
$$ S^{(i)}_t= \mathcal E(\hat{X}^{(i)})_t$$
where $\mathcal E(\cdot )$ is Dolean-Dade exponential. We remark that if $X$ is  a Levy process then $\hat{X}$ is again a Levy process and that
$$d\,S^{(i)}_t=  S^{(i)}_{t-} \,d\hat{X}^{(i)}_t.$$
In addition, for $1\leq i\leq d$ and $t\in[0,T]$
$$ \hat{X}^{(c),i}_t = X^{(c),i}_t $$
$$\nu^ {\hat{X}^{(i)},Q^{*}}= (e^{y_i}-1)\cdot\nu^ {X^{(i)},Q^{*}} .$$
Replacing in Theorem \ref{GR2}   the process $S$ by the process $\hat{X}$ we obtain $Q$-a.s. for all $t\leq T$:
\begin{equation}\label{decmin}
E_{Q^*}[f'(Z^*_T)|\mathcal{F}_t]= x_0 +\sum_{i=1}^d[ \int_0^t\phi ^{(i)}_s S^{(i)}_{s-} d\hat{X}^{(c),Q^{*},i}_s + \int_0^t \int_{\mathbb{R}^{d}}\phi ^{(i)}_s S_{s^-}^{(i)}\,d(\mu^{\hat{X}^{(i)}}-\nu^ {\hat{X}^{(i)},Q^{*}})]
\end{equation}
Then it follows from (\ref{decmin}), Theorem \ref{thdec} and the unicity of decomposition of  martingales on continuous and discontinuous parts, that $Q^{*}-a.s.$, for all $s\leq T$ and all $y\in supp(\nu)$, 
\begin{equation}\label{discont}
H_s(Z^*_{s-},y)=\sum_{i=1}^d \phi^{(i)}_s S^{(i)}_{s-}(e^{y_i}-1)
\end{equation} 
and  for  all $t\leq T$
\begin{equation}\label{cont}
\sum_{i=1}^d \int_0^t \xi_s(Z^*_{s-})\,Z^*_{s-}\,\beta _i^*\,dX^{(c),Q^{*},i}_s = \sum_{i=1}^d \int_0^t \phi^{(i)}_{s} S^{(i)}_{s-}dX^{(c),Q^{*},i}_s.
\end{equation}
We remark that $Q^{*}-a.s.$ for all $s\leq T$
$$H_s( Z^*_s, y) = E_{Q^{*}} (f'( Y^*(y) Z^*_T) - f'( Z^*_T)\,|\, \mathcal F _s).$$
Moreover, $H_s( Z^*_{s-}, y)$ coincide with $H_s( Z^*_s, y)$ in points of continuity of $Z^*$. Taking the sequence  of continuity points of $Z^*$  tending to $T$ and using that $Z_T=Z_{T-}$ ($Q^*$-a.s.) we get that $Q^{*}-a.s.$ for $y\in supp (\nu)$ 
\begin{equation}\label{equality1}
f'(Z^*_{T}Y^*(y))-f'(Z^*_{T})=\sum_{i=1}^d\phi ^{(i)}_{T-} S^{(i)}_{T-}(e^{y_i}-1)
\end{equation}
We fix an arbitrary $y_0\in \stackrel{\circ}{supp}(\nu)$. Differentiating with respect to $y_i$, $i\leq d$, we obtain that 
$$Z^*_{T}\frac{\partial}{\partial y_i}Y^*(y_0)f''(Z^*_{T}Y^*(y_0))=\phi ^{(i)}_{T-}S^{(i)}_{T-}e^{y_{0,i}}$$
We also define : 
$$\Phi(x)=xf''(xY^*(y_0))$$ and
 $$\alpha_i=e^{-y_{0,i}}\frac{\partial}{\partial y_i}Y^*(y_0).$$
We then have $\phi^{(i)}_{T-}S^{(i)}_{T-}=\Phi(Z^*_{T})\alpha_i,$ and inserting this in (\ref{equality1}), we obtain (\ref{eqfmin}).
  
\par 
Taking quadratic variation of the difference of the right-hand side and left-hand side in (\ref{cont}), we obtain that $Q^*-a.s.$ for all $s\leq T$  
$$^{\top}[\xi_s(Z^*_{s-})\,Z^*_{s-}\,\beta^* - S_{s-}\phi_{s}]\,c\,[\xi_s(Z^*_{s-})\,Z^*_{s-}\,\beta^* - S_{s-}\phi_{s}] = 0$$
where by convention $S_{s-}\phi_{s}=(S^{(i)}_{s-}\phi^{(i)}_{s})_{1\leq i \leq d}$. Now, we remark that
$Q^*-a.s.$ for all $s\leq T$
$$Z^*_{s}\, \xi _s( Z^*_s) =  E_{Q^*}(f''( Z^ *_T)Z^*_T \,|\, \mathcal F _s)$$
and that it coincides with $ \xi( Z^*_{s-})$ in continuity points of $Z^*$.
 We take  a set of continuity points  of $Z^*$ which goes to $T$ and we obtain since Levy process has no predictable jumps that $Q^*-a.s.$
$$^{\top}[Z^*_{T}f''(Z^*_{T})\beta^* - S_{T-}\phi_{T-}]\,c\,[Z^*_{T}f''(Z^*_{T})\beta^* - S_{T-}\phi_{T-}] = 0$$
Hence, if $c\neq 0$, 
$$Z^*_{T}f''(Z^*_{T})\beta^*-S_{T-}\phi_{T-}=V$$
where $V\in\mathbb{R}^d$ is a vector which satisfies $cV=0$. Inserting this in (\ref{equality1}) we obtain (\ref{equ}). Theorem \ref{feq} is proved. $\Box$

\noindent \it Proof of Proposition. \ref{support} \rm 
Writing Ito formula we obtain $P$-a.s. for $t\leq T$ :
\begin{equation}
\begin{aligned}
\ln(Z^*_t)=&\sum_{i=1}^d \beta_i^{*}X_t^{(c),i} +\int_0^t \int_{\mathbb{R}^{d}}\ln(Y^*(y))d(\mu^X-\nu^{X,P})
\\&[-\frac{t}{2}{}\,^{\top}\!\beta^*c\beta^* +t\int_{\mathbb{R}^{d}}[\ln(Y^*(y))-(Y^*(y)-1)]\nu(dy)
\end{aligned}
\end{equation}
As we have assumed $Q^*$ to preserve the Levy property, the Girsanov parameters $(\beta^*,Y^*)$ are independent from $(\omega , t)$, and the process  $\ln(Z^*)= (\ln(Z_t^*))_{0\leq t\leq T}$ is a Levy process with the characteristics:\\

$\hspace{1cm}b^{\ln Z^*}= [-\frac{1}{2}{}\,^{\top}\!\beta^*c\beta^* +\int_{\mathbb{R}^{d}}[\ln(Y^*(y))-(Y^*(y)-1)]\,\nu(dy)$,\\

$\hspace{1cm}c^{\ln Z^*} = ^{\top}\!\beta^*c\beta^*$,\\

$\hspace{1cm}d\nu ^{\ln Z^*} = \ln(Y^*(y)\,\nu (dy)$.\\

Now, as soon as ${}^{\top}\beta^* c\beta^*\neq 0$, the continuous component of $\ln(Z^*)$ is non zero, and from  Theorem 24.10 in \cite{Sa} we deduce that $\mbox{supp}(Z_T^*)= \mathbb R ^{+,*}$ and, hence, i).\\
If $Y^*(y)$ is not identically 1 on $\stackrel{\circ}{supp}(\nu)$, then in  (\ref{eqfmin}) the $\alpha _i$, $1\leq i\leq d$, are not all zeros, and hence, the set $supp( \nu^{\ln(Z^*)})= \{\ln Y^*(y), y \in supp(\nu)\}$ contains an interval. It implies that
$\stackrel{\circ}{supp}(\nu ^{\ln Z^*})\neq \emptyset$. Since $0\in supp(\nu)$,  again from (\ref{eqfmin}) it follows that $0\in supp(\nu ^{\ln Z^*})$. Then ii) is a consequence of Theorem 24.10 in \cite{Sa}.$\Box$

\section{So which $f$ can give MEMM preserving Levy property~?}\label{5}
If one considers some simple models, it is not difficult to obtain $f$-divergence minimal equivalent martingale measures for a variety of functions. In particular, one can see that the $f$-divergence minimal measure does not always preserve the Levy property. What can we claim for the functions $f$ such that $f$-divergence minimal martingale measure exists and preserve Levy property?

\begin{theorem}\label{final}
Let $f:\mathbb{R}^{+*}\rightarrow \mathbb{R}$ be a strictly convex function of class $\mathcal{C}^3$ and let $X$ be a Levy process given by its characteristics $(b,c,\nu)$. Assume there exists an $f$-divergence minimal martingale measure $Q^*$ on a time interval $[0,T]$, which preserves the Levy property and belongs to $\mathcal K ^*$.\\
Then, if $supp(\nu)$ is of the non-empty interior,  it contains zero and $Y$ is not identically 1, there exists $a>0$ and $\gamma\in\mathbb{R}$ such that for all $x\in supp(Z^*_T)$, 
$$f''(x)=ax^{\gamma}.$$
If  ${}^{\top}\beta^* c\beta^*\neq 0$ and there exists $y\in supp (\nu )$ such that $ Y^*(y)\neq 1$, then there exist $n\in \mathbb N$, $\gamma\in \mathbb R$, $a>0$ and the real constants $b_i, \tilde{b}_i, 1\leq i\leq n$, such that
$$f''(x)=  ax^{\gamma} + x^{\gamma}\sum_{i=1}^n b_i (\ln(x))^i +\frac{1}{x} \sum_{i=1}^n \tilde{b}_i (\ln(x))^{i-1}$$
\end{theorem}

We deduce this result from the equations obtained in Theorem \ref{feq}. We will successively consider the cases when $\stackrel{\circ}{supp}(\nu)\neq \emptyset$, then when  $c$ is invertible,   and finally when $c$ is not invertible.

\subsection{First case : the interior of $supp(\nu)$ is not empty}
\par \noindent \it Proof of Theorem \ref{final}. \rm We assume that   $\stackrel{\circ}{supp}(\nu)\neq \emptyset$, $0\in supp(\nu)$ , $Y^ *$ is not identically 1 on  ${supp}(\nu)$. According to the  Proposition \ref{support} it implies in both cases 
${}^{\top}\beta^* c\beta^*\neq 0$ and ${}^{\top}\beta^* c\beta^*= 0$, that  $supp (Z_T^*)$ is an interval, say $J$.
Since  the interior of $supp(\nu)$ is not empty, there exist  open non-empty intervals $I_1,...I_d$ such that $I=I_1\times...\times I_d\subseteq \stackrel{\circ}{supp}(\nu)$. Then it follows from Theorem \ref{feq} that for all $(x,y)\in J\times I$, 
\begin{equation}\label{eqdep}
f'(xY^*(y))-f'(x)=\Phi(x)\sum_{i=1}^d \alpha_i(e^{y_i}-1)
\end{equation}
where $\Phi$ is a differentiable on $\stackrel{\circ}{J}$ function and $\alpha\in\mathbb{R}^d$.
 If we now fix $x_0\in \stackrel{\circ}{J}$, we obtain 
$$Y^*(y)=\frac{1}{x_0}(f')^{-1}(f'(x_0)+\Phi(x_0)\sum_{i=1}^d \alpha_i (e^{y_i}-1))$$
and so $Y^*$ is differentiable and monotonous in each variable. Since
$Y^ *$ is not identically 1 on  $\stackrel{\circ}{supp}(\nu)$ we get that $\alpha\neq 0$.
We may now differentiate (\ref{eqdep}) with respect to  $y_i$  corresponding to $\alpha _i\neq 0$, to obtain for all $(x,y)\in J\times I$, 
\begin{equation}\label{neweq}
\Psi(x_0)f''(xY^*(y))=\Psi(x)f''(x_0Y^*(y)),
\end{equation}
where $\Psi(x)=\frac{\Phi(x)}{x}$. Differentiating this new expression with respect to $x$ on the one hand, and with respect to   $y_i$  on the other hand, we obtain the system
\begin{equation}\label{sys}
\begin{cases}&\Psi(x_0)Y^*(y)f'''(xY^*(y))=f''(x_0Y^*(y))\Psi'(x)
\\&\Psi(x_0)xf'''(xY^*(y)) = x_0 f'''(x_0Y^*(y))\Psi(x)
\end{cases}
\end{equation}
In particular, separating the variables, we deduce from this system that there exists $\gamma\in\mathbb{R}$ such that for all $x\in\stackrel{\circ}{J}$,
$$\frac{\Psi'(x)}{\Psi(x)}=\frac{\gamma}{x}.$$ 
Hence, there exists $a>0$ and $\gamma\in\mathbb{R}$ such that for all $x\in \stackrel{\circ}{J}$, $\Psi(x)=ax^{\gamma}$. 
 It then follows from (\ref{neweq}) and (\ref{sys}) that for all $(x,y)\in J\times I$, 
$$\frac{f'''(xY^*(y))}{f''(xY^*(y))}=\frac{\gamma}{xY^*(y)}$$
and hence that $f''(xY^*(y))=a(xY^*(y))^{\gamma}$.
\par We take now the sequence of $(y_m)_{m\geq 1}$, $y_m\in supp(\nu)$, going to zero. Then, the sequence $(Y^*(y_m))_{m\geq 1}$ according to the formula for $Y^*$, is going to 1. Inserting $y_m$ in previous expression and passing to the limit we obtain that for all $x\in \stackrel{\circ}{J}$,
$$\frac{f'''(x)}{f''(x)}=\frac{\gamma}{x}$$
and it proves the result on $\stackrel{\circ}{supp}(Z_T^*)$. The final result on $supp( Z^*_T)$ can be proved again by limit passage. $\Box$

\subsection{Second case : $c$ is invertible and  $\nu$ is nowhere dense}
In the first case, the proof relied on differentiating the function $Y^*$. This is of course no longer possible when the support of $\nu$ is nowhere dense. Howerever, since  ${}^{\top}\beta^* c\beta^*\neq 0$, we get from Proposition \ref{support} that $supp(Z^*)= \mathbb R^ {+,*}$.
Again from Theorem \ref{feq}  we have for all $x>0$ and  $y\in supp(\nu)$, 
\begin{equation}\label{eqdep2}
f'(xY^*(y))-f'(x)=xf''(x)\sum_{i=1}^d \beta _i^{*}(e^{y_i}-1).
\end{equation}
We will distinguish two similar cases: $b>1$ and $0<b<1$. For $b>1$
we  fix $\epsilon$, $0<\epsilon<1$, and we introduce for $a\in\mathbb{R}$ the following vector space:
$$V_{a,b}=\{\phi\in \mathcal{C}^1([\epsilon(1\wedge b),\frac{1\vee b}{\epsilon}]), \text{ such that for }x\in[\epsilon,\frac{1}{\epsilon}], \phi(bx)-\phi(x)=ax\phi'(x)\}
$$
with the norm $$|\!|\phi |\!|_{\infty}= \sup_{x\in[\epsilon,\frac{1}{\epsilon}]}|\phi (x)| + \sup_{x\in[\epsilon,\frac{1}{\epsilon}]} | \phi(bx)|$$
It follows from (\ref{eqdep2})  that $f'\in V_{a,b}$ with $b= Y^*(y)$ and $a= \sum_{i=1}^d \beta _i^{*}(e^{y_i}-1)$. The condition that there exist $y\in supp (\nu)$ such that $Y^*(y)\neq 1$ insure that $\sum_{i=1}^d \beta _i^{*}(e^{y_i}-1)\neq 0$.

\begin{lemma}
If $a\neq 0$ then $V_{a,b}$ is a finite dimensional closed in $|\!|\cdot |\!|_{\infty}$ vector space.
\end{lemma}
\noindent \it Proof \rm It is easy to verify that $V_{a,b}$ is a vector space. We show that $V_{a,b}$ is a closed vector space : if we consider a sequence $(\phi_n)_{n\geq 1}$ of elements of $V_{a,b}$ which converges to a function $\phi$, we denote by $\psi$ the function such that $\psi(x)=\frac{\phi(bx)-\phi(x)}{ax}$. We then have
$$\lim_{n\to+\infty}||\phi'_n-\psi||_{\infty}\leq \frac{1}{\epsilon |a|(1\wedge b)}\lim_{n\to+\infty}||\phi_n-\phi||_{\infty}=0$$
Therefore, $\phi$ is differentiable and we have $\phi'=\psi $. Therefore, $\phi$ is of class $\mathcal{C}^1$ and belongs to $V_{a,b}$. Hence, $V_{a,b}$ is a closed in $|\!|\cdot |\!|_{\infty}$ vector space. Now, for $\phi\in V_{a,b}$ and $x,y\in[\epsilon ,\frac{1}{\epsilon}]$, we have 
$$|\phi(x)-\phi(y)|\leq \sup_{u\in[\epsilon,\frac{1}{\epsilon}]}|\phi'(u)||x-y|\leq \sup_{u\in[\epsilon,\frac{1}{\epsilon}]}\frac{|\phi(bu)-\phi(u)|}{|au|}|x-y|\leq \frac{||\phi||_{\infty}}{|a|\epsilon}|x-y|$$
Therefore, the unit ball of $V_{a,b}$ is equi-continuous, hence, by Ascoli theorem, it is relatively compact, and now it follows from the Riesz Theorem that $V_{a,b}$ is a finite dimensional vector space. $\Box$

\vskip 0.2cm  We now show that elements of $V_{a,b}$ belong to a specific class of functions. 
\begin{lemma}
All elements of $V_{a,b}$ are solutions to a Euler type differential equation, that is to say there exists $m\in\mathbb{N}$ and  real numbers $(\rho_i)_{0\leq i\leq m}$ such that 
\begin{equation}\label{span}
\sum_{i=0}^m \rho_i x^i \phi^{(i)}(x)=0.
\end{equation}
\end{lemma}
\noindent \it Proof \rm It is easy to see from the definition of $V_{a,b}$ that if $\phi\in V_{a,b}$, then the function $x\mapsto x\phi'(x)$ also belongs to $V_{a,b}$. If we now denote by $\phi^{(i)}$ the derivative of order $i$ of $\phi$, we see that the span of $(x^i\phi^{(i)}(x))_{i\geq 0}$ must be a subvector space of $V_{a,b}$ and in particular a finite dimensional vector space. In particular, there exists $m\in\mathbb{N}$ and real constants $(\rho_i)_{0\leq i\leq m}$ such that (\ref{span}) holds. $\Box$

\vskip 0.2cm \noindent \it Proof of Theorem \ref{final} \rm The previous result applies in particular to the function $f'$ since $f'$ verify (\ref{eqdep2}). As a consequence, $f'$ satisfy Euler type differential equation. It is known that the change of variable $x=\exp (u)$ reduces this equation to a homogeneous differential equation of order $m$ with constant coefficients. It is also known that the solution of such equation can be written as a linear combination of the solutions corresponding to different roots of characteristic polynomial. These solutions being linearly independent, we need only to considerer a generic one, say $f'_{\lambda}$, $\lambda$ being the root of characteristic polynomial. If the root of characteristic polynomial $\lambda$ is real and of the multiplicity $n$, $n\leq m$, then
$$f'_{\lambda}(x)=  a_0 x^{\lambda}+x^{\lambda} \sum_{i=1}^n b_i (\ln(x))^i$$
and if this root is complex then
$$f'_{\lambda}(x)=x^{Re(\lambda)}\sum_{i=0}^n [c_i \cos(\ln(Im(\lambda)x))+d_i \sin(\ln(Im (\lambda)x))]\ln(x)^i$$
where  $a_0,b_i, c_i, d_i$ are real constants.
Since $f'$ is increasing,  we must have for all $i\leq n$, $c_i=d_i=0$. But $f$ is strictly convex and the last case is excluded. Putting
$$f'_{\lambda}(x)= a_0 x^{\lambda}+x^{\lambda} \sum_{i=1}^n b_i (\ln(x))^i$$
into the  equation
\begin{equation}\label{fi}
f'(bx)-f'(x)= a x f''(x)
\end{equation}
we get using linear independence of mentioned functions that

\begin{equation}\label{m0}
a_0(b^{\lambda} - a\lambda -1) +b^{\lambda}\sum_{i=1}^n b_i(\ln b)^i - ab_1=0 
\end{equation}
and that for all $1\leq i\leq n$,
\begin{equation}\label{m}
\sum_{k=i}^n b^{\lambda} b_k C_k^{i}(\ln(b))^{k-i}-b_i(1+a \lambda) - a b_{i+1}(i+1)=0 
\end{equation}
with $b_{n+1}=0$.
We remark that 
the matrix corresponding to (\ref{m}) is   triangular matrix $M$ with $b^{\lambda}-1-a \lambda $ on the diagonal. If $b^{\lambda}-1-a \lambda \neq 0$, then  the system of equations has unique solution. This solution should also verify:for all $x>0$
\begin{equation}\label{m1}
f''_{\lambda}(x)>0
\end{equation}
If $b^{\lambda}-1-a \lambda =0$ , then $rang(M)=0$, and
$b_i$ are free constants. Finally, we conclude that there exist  a solution
$$f'_{\lambda}(x)= ax^{\lambda} + x^{\lambda}\sum_{i=1}^n b_i (\ln(x))^i$$
verifying (\ref{m1}) with any $\lambda$ verifying $b^{\lambda}-1-a \lambda=0$.
$\Box$

\subsection{Third case : $c$ is non invertible  and  $\nu$ is nowhere dense}
We finally consider the case of Levy models which have a continuous component but for which the matrix $c$ is not invertible. It follows from Theorem \ref{feq} that in this case we have for all $x\in supp(Z^*)$ and $y\in supp(\nu)$
\begin{equation}\label{third} 
f'(xY^*(y))-f'(x)=xf''(x)\sum_{i=1}^d \beta _i^{*}(e^{y_i}-1)-\sum_{j=1}^d V_j(e^{y_j}-1)
\end{equation}
where $cV=0$. 
\vskip 0.2cm \noindent \it Proof of Theorem \ref{final} \rm First of all, we note that if $f'$ satisfies (\ref{third}) then $\phi:x\mapsto xf''(x)$ satisfies (\ref{fi}). The conclusions of the previous section then hold for $\phi$. $\Box$

\section{Minimal equivalent measures when $f''(x)=ax^{\gamma}$}\label{6}
Our aim in this section is to consider in more detail the class of minimal martingale measures for the functions which satisfy $f''(x)=ax^{\gamma}$. First of all, we note that these functions are those for which there exists $A>0$ and real $B$,$C$ such that 
$$f(x)=Af_{\gamma}(x)+Bx+ C$$
where
\begin{equation}\label{formoff}
f_{\gamma}(x)=\begin{cases}& c_{\gamma} x^{\gamma+2} \text{ if }\gamma\neq -1,-2,
\\& x\ln(x) \text{ if }\gamma=-1,
\\&-\ln(x) \text{ if }\gamma=-2.
\end{cases}
\end{equation}
and $c_{\gamma}= \mbox{sign}[(\gamma +1)/(\gamma+2)]$. 
In particular, the minimal measure for $f$ will be the same as that for $f_{\gamma}$. Minimal measures for the different functions $f_{\gamma}$ have been well studied. It has  been shown in \cite{K1}, \cite{ES}, \cite{JKM} that in all these cases, the minimal measure, when it exists, preserves the Levy property.\\ 
 Sufficient conditions for the existence of a minimal measure and an explicit expression of the associated Girsanov parameters have been given in the case of relative entropy in \cite{FM},\cite{HS} and for power functions in \cite{JKM}. It was also shown in \cite{HS} that these conditions are in fact necessary in the case of relative entropy or for power functions when $d=1$. Our aim in this section is to give a unified expression of such conditions for all functions which satisfy $f''(x)=ax^{\gamma}$ and to show that, under some conditions, they are necessary and sufficient, for all $d$-dimensional Levy models.\\
 We have already mentioned that $f$-divergence minimal martingale measures play an important role in the determination of utility maximising strategies. In this context, it is useful to have further invariance properties for the minimal measures such as scaling and time invariance properties. This is  the case when $f''(x)=ax^{\gamma}$. \\

\begin{theorem}\label{exist}
Consider a Levy process $X$ with characteristics $(b,c,\nu)$  and let $f$ be a function such that $f''(x)=ax^{\gamma}$, where $a>0$ and $\gamma\in\mathbb{R}$. Suppose that $c\neq 0$ or $\stackrel{\circ}{supp}(\nu)\neq \emptyset$. Then there exists an $f$-divergence minimal equivalent  to $P$ martingale measure $Q$ preserving Levy properties if and only if there exist $\gamma , \beta\in\mathbb{R}^d$ and measurable function  $Y : \mathbb R ^ d\setminus \{0\}\rightarrow \mathbb R^ {+}$  such that
\begin{equation}\label{Y}
Y(y)=(f')^{-1}(f'(1)+\sum_{i=1}^d \gamma _i(e^{y_i}-1))
\end{equation}
and such that the following properties hold:
\begin{equation}\label{cdsec1}
Y(y)> 0 \,\,\,\nu-a.e.,
\end{equation}
\begin{equation}\label{cdsec2}
\sum_{i=1}^d \int_{|y|\geq 1}(e^{y_i}-1)Y(y)\nu(dy)<+\infty.
\end{equation}
\begin{equation} \label{cdsec3}
b+\frac{1}{2}diag(c)+c\beta+\int_{\mathbb{R}^{d}}((e^y-1)Y(y)-h(y))\nu(dy)=0.
\end{equation}
If such a measure exists the Girsanov parameters associated with $Q$ are $\beta$ and $Y$, and this measure is scale and time invariant.
\end{theorem}
We begin with some technical lemmas.\\

\begin{lemma}\label{measures}
Let $Q$ be the measure preserving Levy property. Then, $Q_T\sim P_T$ for all $T>0$ iff
\begin{equation}
Y(y)> 0 \,\,\,\nu-a.e.,
\end{equation}
\begin{equation}\label{hellinger}
 \int_{\mathbb{R}^{d}}(\sqrt{Y(y)}-1)^2\nu(dy)<+\infty .
\end{equation}
\end{lemma}

\it Proof \rm \, See Theorem 2.1, p. 209 of \cite{JSh}.$\Box$\\

\begin{lemma}\label{integrability}
Let $Z_T = \frac{d Q_T}{d P_T}$. Under $Q_T\sim P_T$, the condition $E_P | f(Z_T)| < \infty$ is equivalent to
\begin{equation}\label{predictable}
 \int_{\mathbb{R}^{d}} [f(Y(y))-f(1)-f'(1)(Y(y)-1)] \nu(dy)<+\infty 
\end{equation}
\end{lemma}

\it Proof\rm\, In our particular case, $E_P | f(Z_T)| < \infty$ is equivalent to the existence of $E_P  f(Z_T)$. 
We use Ito formula to express this integrability condition in predictable terms. Taking for $n\geq 1$ stopping times 
$$ s_n= \inf\{ t\geq 0 : Z_t> n \,\mbox{or}\, Z_t < 1/n\}$$
where $\inf\{\emptyset\}= +\infty$,
we get for $\gamma\neq -1,-2$ and $\alpha = \gamma +2$ that $P$-a.s.
$$Z_{T\wedge s_n}^{\alpha} = 1 + \int _0^{T\wedge s_n}\alpha\,Z_{s-}^{\alpha}\beta d X_s^c 
 + \int _0^{T\wedge s_n}\int_{\mathbb R ^ d} Z_{s-}^{\alpha}(Y^{\alpha}(y) -1) (\mu^X - \nu^ {X,P})(ds , dy) $$
$$+ \int _0^{T\wedge s_n}\int_{\mathbb R ^ d} Z_{s-}^{\alpha}[Y^{\alpha}(y) -1 - \alpha(Y(y) -1)] d s\,\, \nu (dy)$$
Hence,
\begin{equation}\label{formula}
Z_{T\wedge s_n}^{\alpha} = \mathcal E( N^{(\alpha )} + A^{(\alpha )})_{T\wedge s_n}
\end{equation}
where
$$N^{(\alpha )}_t=  \int _0^t\alpha\,\beta d X_s^c +
\int _0^t(Y^{\alpha}(y) -1)  (\mu^X - \nu^ {X,P})(ds , dy)$$
and
$$A^{(\alpha )}_t= \int _0^t\int_{\mathbb R ^ d}[Y^{\alpha}(y) -1 - \alpha(Y(y) -1)] d s\,\,  \nu (dy)$$
Since $[N^{(\alpha )}, A^{(\alpha )}]_t =0$ for each $t\geq 0$ we have
$$Z_{T\wedge s_n}^{\alpha} = \mathcal E( N^{(\alpha )})_{T\wedge s_n}\mathcal E ( A^{(\alpha )})_{T\wedge s_n}$$
If $E_P Z_T^{\alpha}< \infty$, then by Jensen inequality
$$0 \leq Z_{T\wedge s_n}^{\alpha} \leq E_P (Z_{T}^{\alpha}\,|\,\mathcal F _{T\wedge s_n})$$
and since the right-hand side of this inequality form uniformly integrable sequence, $(Z_{T\wedge s_n}^{\alpha})_{n\geq 1}$
is also uniformly integrable.
We remark that in the case $\alpha >1$ and $\alpha <0$, $A^{(\alpha )}_t\geq 0$ for all $t\geq 0$ and
$$\mathcal E ( A^{(\alpha )})_{T\wedge s_n} = \exp ( A^{(\alpha )}_{T\wedge s_n})\geq 1.$$
It means that $( \mathcal E( N^{(\alpha )})_{T\wedge s_n})_{n\in\mathbb N^*}$ is uniformly integrable
and 
\begin{equation}\label{du1}
E_P(Z_T^{\alpha}) = \exp (A^{(\alpha )}_T)
\end{equation}
If (\ref{predictable}) holds, then by Fatou lemma and since $\mathcal E( N^{(\alpha )})$ is a local martingale
we get
$$E_P(Z_T^{\alpha}) \leq \underline{\lim}_{n\rightarrow \infty}E_P Z_{T\wedge s_n}\leq \exp (A^{(\alpha )}_T)$$
For $0<\alpha <1$, we have again
$$Z_{T\wedge s_n}^{\alpha} = \mathcal E(N^{(\alpha )})_{T\wedge s_n}\mathcal E (A^{(\alpha )})_{T\wedge s_n}$$
with uniformly integrable sequence $(Z_{T\wedge s_n}^{\alpha})_{n\geq 1}$.
Since
$$\mathcal E ( A^{(\alpha )})_{T\wedge s_n} = \exp ( A^{(\alpha )}_{T\wedge s_n})\geq \exp ( A^{(\alpha )}_{T}),$$
the sequence
 $( \mathcal E(N^{(\alpha )})_{T\wedge s_n})_{n\in\mathbb N^*}$ is uniformly integrable and
\begin{equation}
E_P(Z_T^{\alpha}) = \exp (A^{(\alpha )}_T).
\end{equation} 
For $\gamma = -2$ we have that $f(x) = x\ln (x)$ up to linear term and 
$$Z_{T\wedge s_n}\ln (Z_{T\wedge s_n}) = $$
$$\int _0^{T\wedge s_n} (\ln (Z_{s-}) +1)Z_{s-}\beta d X_s^c
 +
\int _0^{T\wedge s_n} \int _{\mathbb R ^ d}[\ln (Z_{s-}Y(y)) -  \ln (Z_{s-})] (\mu^X - \nu^ {X,P})(ds , dy)$$
$$ + \int _0^{T\wedge s_n}\int _{\mathbb R ^ d} Z_{s-}[Y(y)\ln (Y(y)) -Y(y) +1] ds \,\nu (dy)$$
Taking mathematical expectation we obtain:
\begin{equation}\label{log}
E_P [Z_{T\wedge s_n}\ln (Z_{T\wedge s_n})] = E_P \int _0^{T\wedge s_n}\int _{\mathbb R ^ d} Z_{s-}[Y(y)\ln (Y(y)) -Y(y) +1] d s \,\, \nu (dy)
\end{equation}
If $E_P [Z_{T}\ln (Z_{T})]<\infty$, then the sequence $( Z_{T\wedge s_n}\ln (Z_{T\wedge s_n}))_{n\in\mathbb N^*}$ is uniformly  integrable
 and $E_P(Z_{s-})=1$ we obtain applying Lebesgue convergence theorem that
\begin{equation}\label{du2}
E_P [Z_{T}\ln (Z_{T})] = T\int _{\mathbb R ^ d} [Y(y)\ln (Y(y)) -Y(y) +1] \nu (dy)
\end{equation}
and this implies (\ref{predictable}). If (\ref{predictable}), then by Fatou lemma from (\ref{log}) we deduce that
$E_P [Z_{T}\ln (Z_{T})]<\infty$.\\
For $\gamma =-1$, we have $f(x)= -\ln (x)$ and exchanging $P$ and $Q$ we get:
$$E_P[-\ln (Z_T)]= E_Q [\tilde{Z}_{T}\ln (\tilde{Z}_{T})] = T\int _{\mathbb R ^ d} [\tilde{Y}(y)\ln (\tilde{Y}(y)) -\tilde{Y}(y) +1] \nu^Q (dy)$$ 
where  $\tilde{Z}_T= 1/Z_T$ and $\tilde{Y}(y)= 1/Y(y)$.
 But $\nu ^Q(dy) = Y(y) \nu(dy)$ and, finally,
\begin{equation}\label{du3}
E_P [-\ln (Z_{T})] = T\int _{\mathbb R ^ d} [-\ln (Y(y)) +Y(y) -1] \nu (dy)
\end{equation}  
which implies (\ref{predictable}).$\Box$\\

\begin{lemma}\label{equivalence} If the second Girsanov parameter $Y$ has a particular form (\ref{Y})
then the condition
\begin{equation}\label{cdsec2b}
\sum_{i=1}^d \int_{|y|\geq 1}(e^{y_i}-1)Y(y)\nu(dy)<+\infty
\end{equation}
implies the conditions (\ref{hellinger}) and (\ref{predictable}).
\end{lemma}
\it Proof\,\rm \,\, We can cut each  integral in (\ref{hellinger}) and (\ref{predictable}) on two parts and integrate on the sets $\{|y|\leq 1\} $ and $\{|y|>1\}$. Then we can use a particular form of $Y$ and conclude easily writing Taylor expansion of order 2. $\Box$\\

\noindent \it Proof of Theorem \ref{exist}\, Necessity \rm  We suppose that there exist $f$-divergence minimal equivalent martingale measure $Q$ preserving Levy property of $X$. Then, since  $Q_T\sim P_T$, the conditions (\ref{cdsec1}),
(\ref{hellinger}) follow from Theorem 2.1, p. 209 of \cite{JSh}. From Theorem 3 we deduce that (\ref{Y}) holds. Then, the condition (\ref{cdsec2}) follows from the fact that $S$ is a martingale under $Q$. Finally, the condition (\ref{cdsec3}) follows from Girsanov theorem since $Q$ is a martingale measure and, hence,  the drift of $S$ under $Q$ is zero.\\ 

\it Sufficiency \rm We take $\beta$ and  $Y$ verifying the conditions (\ref{cdsec1}),(\ref{cdsec2}),(\ref{cdsec3}) and we construct
\begin{equation}\label{mart1}
M_t= \sum_{i=1}^d \int_0^t \beta^{(i)} dX^{c, (i)}_s+\int_0^t \int_{\mathbb{R}^{d}}(Y(y)-1)(\mu^{X}-\nu^{X, P})(ds, dy)
\end{equation}
As known from Theorem 1.33, p.72-73,  of \cite{JSh}, the last stochastic integral is well defined if
$$C(W)= T  \int_{\mathbb{R}^{d}}(Y(y)-1)^2I_{\{|Y(y)-1|\leq 1\}} \nu (dy) < \infty ,$$
$$C(W')= T  \int_{\mathbb{R}^{d}}|Y(y)-1|I_{\{|Y(y)-1|> 1\}} \nu (dy) < \infty .$$
But the condition (\ref{cdsec2}), the relation (\ref{Y}) and  lemma \ref{equivalence} implies (\ref{hellinger}). So, $(Y-1)\in G_{loc}(\mu ^ X)$ and
$M$ is local martingale. Then we take
$$Z_T= \mathcal E( M)_T$$
and this defines the measure $Q_T$ by its Radon-Nikodym density.  Now, the conditions (\ref{cdsec1}),(\ref{cdsec2}) together with the relation (\ref{Y}) and  lemma \ref{equivalence} imply (\ref{hellinger}), and, hence,
from lemma \ref{measures} we deduce $P_T\sim Q_T$.\\
We show that $E_P|f(Z_T)|<\infty $. Since $P_T\sim Q_T$, the lemma \ref{integrability} gives needed integrability condition. \\
Now, since (\ref{cdsec3}) holds, $Q$ is martingale measure, and it remains to show that $Q$ is indeed
$f$-divergence minimal. For that we take any equivalent martingale measure $\bar{Q}$ and we show that
\begin{equation}\label{ineq}
E_Q f'(Z_T) \leq E_{\bar{Q}} f'(Z_T).
\end{equation}
If the mentioned  inequality holds, the Theorem \ref{GR1}
implies that $Q$ is a minimal.\\
In the case $\gamma \neq -1,-2$ we obtain from (\ref{formula}) replacing $\alpha$ by $\gamma +1$:
$$Z_T^{\gamma +1} = \mathcal E(N^{(\gamma +1)})_T \,\,\exp(A_T^{(\gamma +1)})$$
and using a particular form of $f'$ and $Y$  we get that for $0\leq t\leq T$
$$N_t^{(\gamma +1)} =\sum_{i=1}^d \theta ^{(i)} \hat{X}^{(i)}_t$$ 
where $ \theta = \beta$ if $c\neq 0$ and $\theta = \gamma$ if $c=0$, and $\hat{X}^{(i)}$ is a stochastic logarithm of $S^{(i)}$. So, $\mathcal  E(N^{(\gamma +1)})$ is a local martingale and we get
$$E_{\bar{Q}} Z_T^{\gamma +1} \leq \exp(A_T^{(\gamma +1)}) = E_Q Z_T^{\gamma +1}$$
and, hence, (\ref{ineq}).\\
In the case $\gamma = -1$
we prove using again a particular form of $f'$ and $Y$ that
$$f'(Z_T) = E_Q( f'(Z_T)) + \sum_{i=1}^d {\theta}^{(i)} \hat{X}^{(i)}_T$$ 
with $ \theta = \beta$ if $c\neq 0$ and $\theta = \gamma$ if $c=0$.
Since $E_{\bar{Q}} \hat{X}_T=0$ we get that
$$E_{\bar{Q}}( f'(Z_T)) \leq E_Q( f'(Z_T))$$
and it proves that $Q$ is $f$-divergence minimal.\\
The case $\gamma =-2$ can be considered in similar way.\\

Finally, note that the conditions which appear in Theorem \ref{exist} do not depend in any way on the time interval which is considered and, hence, the minimal measure is time invariant. Furthermore,  if $Q^*$ is $f$-divergence minimal, the equality 
$$f(cx)= Af(x) + B x+ C$$
with $A,B,C$ constants, $A>0$,
gives
$$E_P[f(c\frac{d\bar{Q}}{dP})]=A E_P[f(\frac{d\bar{Q}}{dP})]+B +C \geq A E_P[f(\frac{dQ}{dP})]+ B + C=E_P[f(c\frac{dQ}{dP})]$$
and $Q$ is scale invariant. $\Box$

\subsection{ Example}
We now give an example of a Levy model and a convex function which does not satisfy $f''(x)=ax^{\gamma}$ yet preserves the Levy property. 
 We consider the function $f(x)=\frac{x^2}{2}+x\ln(x)-x$ and the $\mathbb{R}^2$-valued Levy process given by 
$X_t=(W_t+\ln(2)P_t,W_t+\ln(3)P_t-t)$, where $W$ is a standard one-dimensional Brownian motion and $P$ is a  standard one-dimensional Poisson process. Note that the covariance matrix 
$$c= \left(
\begin{array}{ll}
1&1\\
1&1
\end{array}\right)
$$
 is not invertible. The support of the Levy measure is the singleton $a=(\ln(2),\ln(3))$, and is in particular nowhere dense. Let $Q$ be a martingale measure for this model, and $(\beta,Y)$ its Girsanov parameters, where ${}^{\top}\beta=(\beta_1,\beta_2)$. In order for $Q$ to be a martingale measure preserving Levy property, we must have 
\begin{equation}\label{eqexgirs}
\begin{aligned}
&\ln(2)+\frac{1}{2}+\beta_1+\beta_2+Y(a)=0,
\\& \ln(3)-\frac{1}{2}+\beta_1+\beta_2+2Y(a)=0,
\end{aligned}
\end{equation}
and, hence, $Y(a)=1-\ln(\frac{3}{2})$.
Now, it is not difficult to verify using Ito formula that the measure $Q$ satisfy:
$E_P Z_T^ 2 <\infty$ and, hence, $E_P |f(Z_T)|<+\infty$ . Moreover, the conditions (\ref{cdsec1}),(\ref{cdsec2})
are satisfied meaning that $P_T\sim Q_T$.\\
Furthermore, in order for $Q$ to be minimal  we must have according to  Theorem \ref{feq}:
$$f'(xY(y)) - f'(x)= x f''(x)\sum_{i=1}^{2}\beta_i(e^{a_i}-1)+\sum_{i=1}^{2} v_i(e^{a_i}-1)$$
with $a_1=\ln 2, a_2=\ln 3$ and $V= \,^{\top}(v_1,v_2)$ such that $cV=0$. We remark that $v_2=-v_1$.
Then for $x\in supp(Z_T)$
$$\ln(Y(a))+x(Y(a)-1)=(x+1)(\beta_1+2\beta_2)-v_1$$
and since $supp(Z_T)=\mathbb R^{+,*}$ we must have 
$$\beta_1+2\beta_2=Y(a)-1 \text{ and }\beta_1+2\beta_2-v_1=\ln(Y(a))$$
Using (\ref{eqexgirs}), this leads to 
$$\begin{cases}
&v_1=-\ln(1-\ln(\frac{3}{2}))-\ln(\frac{3}{2})
\\&\beta_1=3\ln(3)-5\ln(2)-3
\\&\beta_2=\frac{3}{2}+3\ln(2)-2\ln(3)
\end{cases}$$
We now need to check that the martingale measure given by these Girsanov parameters is indeed minimal. Note that the decomposition of Theorem \ref{thdec} can now be written
$$f'(Z_T)=E_Q[f'(Z_T)]+\sum_{i=1}^2 \int_0^T \left[\beta_i ( \frac{1}{Z_{s-}} +E_Q[Z_{T-s}])+v_i\right] \frac{dS^i_s}{S^i_{s-}}$$
But  for $s\geq 0$
$$  \frac{dS^i_s}{S^i_{s-}}= \hat{X}_s^i$$
and right-hand side of previous equality is a local martingale with respect to any martingale measure
$\bar{Q}$. Taking a localising sequence and then the expectation with respect to $\bar{Q}$ we get after limit passage that
 $$E_{\bar{Q}}[f'(Z_T)]\leq E_{Q}[f'(Z_T)],$$ and so, it follows from Theorem \ref{GR1} that the measure $Q$ is indeed minimal. 
 
 \vskip 0.2cm
\noindent \section{Acknowledgements} \rm  This work is supported in part by ECOS project M07M01  and by ANR-09-BLAN-0084-01 of the Department of Mathematics of Angers's University.


\end{document}